\documentclass[12pt]{article}
\usepackage[utf8]{inputenc}
\usepackage{amsmath}
\usepackage{amssymb}
\usepackage{amsthm}

\title{A criterion for $p$-henselianity in characteristic $p$}

\author{Zo\'e Chatzidakis\thanks{Partially supported by ANR-13-BS01-0006} and Milan Perera
}

\date{}

\begin{document}
\maketitle

\begin{abstract}
Let $p$ be a prime. In this paper we give a proof of the following
result: A valued field $(K,v)$ of characteristic $p > 0$ is
$p$-henselian if and only if every element of strictly positive
valuation if of the form $x^p - x$ for some $x \in K$. 
\end{abstract}


\section*{Preliminaries}

Throughout this paper, all fields have characteristic $p > 0$. First we recall some definitions and notations. Let $\mathcal{O}_v := \{x \in K \mid v(x) \geq 0\}$ be the valuation ring associated with $v$. It is a local ring, and $\mathcal{M}_{v} := \{x \in K \mid v(x) > 0 \}$ is its maximal ideal. Let $\overline{K}_v := \mathcal{O}_v/\mathcal{M}_v = \{\overline{a} = a + \mathcal{M}_v \mid a \in \mathcal{O}_v\}$ be the residue field (or $\overline{K}$ when there is no danger of confusion). We let $K(p)$ denote the \emph{compositum} of all finite Galois extensions of $K$ of degree a power of $p$.

A valued field $(K,v)$ is $p$-henselian if $v$ extends uniquely to
$K(p)$. Equivalently (see \cite{EP}, Thm 4.3.2), if it satisfies a
restricted version of Hensel's lemma (which we call $p$-Hensel lemma) :
$K$ is $p$-henselian if and only if every polynomial $P \in
\mathcal{O}_v[X]$ which splits in $K(p)$ and with residual image in
$\overline{K}_v[X]$ having a simple root $\alpha$ in $\overline{K}_v$,
has a root $a$ in $\mathcal{O}_v$ with $\overline{a} =
\alpha$. Furthermore, another result (see \cite{EP}, Thm 4.2.2) shows
that $(K,v)$ is $p$-henselian if and only if $v$ extends uniquely to
every Galois extension of degree $p$.\\[0.1in] 
The aim of this note is to give a complete proof of the following
result:\\[0.1in]  
{\bf Theorem}. Let $(K,v)$ be a valued field. $(K,v)$ is $p$-henselian \emph{if and only if} $ \mathcal{M}_v \subseteq \{ x^p-x \mid x \in K \}$.\\[0.1in]
This result was announced in \cite{Ko}, Proposition~1.4, however the proof was not complete. The notion of $p$-henselianity is important in the study of fields with definable valuations, and in particular it is important to show that the property of $p$-henselianity is an elementary property of valued fields.

The proof we give is elementary, and uses extensively pseudo-convergent
sequences and their properties. Recall that a sequence $\{a_\rho\}_{\rho < \kappa} \in K^{\kappa}$ indexed by an ordinal $\kappa$ is said to be \emph{pseudo-convergent} if for all $\alpha < \beta < \gamma < \kappa$ : \begin{equation}
v(a_\beta - a_\alpha) < v(a_\gamma - a_\beta).
\end{equation}

A pseudo-convergent sequence $\{a_\rho\}_{\rho < \kappa}$ is called \emph{algebraic} if there is a polynomial $P$ in $K[X]$ such that $v\big(P(a_\alpha)\big) <  v\big(P(a_\beta)\big)$ ultimately for all $\alpha < \beta$, \emph{i.e}: \begin{equation}
\exists \lambda < \kappa \forall \alpha, \beta < \kappa \quad (\lambda < \alpha < \beta) \Rightarrow v\big(P(a_\alpha)\big) <  v\big(P(a_\beta)\big).
\label{alg}
\end{equation} Otherwise, it is called {\em transcendental}.

We assume familiarity with the properties of pseudo-convergent
sequences, see \cite{Ka} for more details, and in particular Theorem 3,
Lemmas 4 and 8. 

\section*{Proof of the theorem}

First, we prove a lemma in order to restrict our study to immediate
extensions: \\[0.1in]
{\bf Observation}. {\em  
Let $(K,v)$ be a valued field and $(L,w)$ be a Galois extension of degree a prime $\ell$. Then, if $(L, w) / (K,v)$ is residual or ramified, $w$ is the unique extension of $v$ to $L$.}

\begin{proof}
The fundamental equality of valuation theory (see \cite{EP}, Thm 3.3.3) tells us that if $L$ is a Galois extension of $K$, then
\begin{equation}
[L:K] = e(L/K)f(L/K)gd
\end{equation} where $e(L/K)$ is the ramification index, $f(L/K)$ the residue index, $g$ the number of extensions of $v$ to $L$ and $d$, the defect, is a power of $p$.

Thus, as $\ell$ is a prime, if $e(L/K)f(L/K) > 1$, then necessarily $g = d = 1$, and in particular, $v$ has a unique extension to $L$.
\end{proof}

Now, let us prove the result announced in the preliminaries: \\[0.1in]
%
{\bf Theorem}. {\em Let $(K, \mathcal{O}_{v})$ be a valued field of characteristic $p$. Then, $(K, \mathcal{O}_{v})$ is $p$-henselian \emph{if and only if} $\mathcal{M}_{v} \subseteq K^{(p)} := \{x^{p} - x \mid x \in K\}$.}

\begin{proof}
The forward direction is an immediate application of the $p$-Hensel Lemma.\newline
Conversely, assume that $\mathcal{M}_{v} \subseteq K^{(p)} := \{x^{p} -
x \mid x \in K\}$. Every Galois extension of $K$ of degree $p$ is an
Artin-Schreier extension, \emph{i.e} is generated over $K$ by a root $a$
of a polynomial $X^p - X - b = 0$, with $b \in K \setminus K^{(p)}$. The
previous observation gives us the result when $K(a)/K$ is not
immediate. Let $L$ be an immediate Galois extension of degree $p$ and
$\tilde{v}$ an extension of $v$ to $L$ (hence with the same value group
$\Gamma$ and residue field $\overline{L} = \overline{K}$ as $K$). We can
write $L = K(a)$ where $a^p - a = b \in K \setminus K^{(p)}$.\\[0.1in]
Step 1: (Claim) The set $C = \{v(x^p - x - b) \mid x \in K
\} = v\big(K^{(p)} - b\big)$ is contained in $\Gamma_{<0}$ 
and has no last element.

First observe that $C \subseteq \Gamma_{\leq 0}$ : if $v(c^p - c - b) > 0$, then the equation $X^p - X + (c^p - c - b)$ has a root in $K$, so that $(a-c) \in K$: 
contradiction. Let $\gamma \in \Gamma$, $d \in K$ such 
that $v(d^p - d - b) = \gamma$. As $L/K$ is immediate 
there is $c \in K$ such that $\tilde{v}\big(a - (d+c)\big) 
> \tilde{v}(a-d)$. If $\tilde{v}(a-d) = 0$ then $\tilde{v}
\big(a - (d + c) \big) > 0$ and $\big((d+c)^p - (d+c) - b 
\big) = (d+c-a)^p - (d+c-a)$ in $\mathcal{M}_v$, which as 
above give a contradiction. Hence $\tilde{v}(a-d) < 0$, and from $d^p - d - b = (d - a)^p - (d - a)$, we deduce 
that $\gamma = p\tilde{v}(a-d) < 0$, and $v\big((d+c)^p - 
(d+c) - b \big) = p\Big(\tilde{v}\big(a-(d+c)\big)\Big) > 
\gamma$. This shows the claim.\\[0.1in]
Step 2: We extract a strictly well-ordered increasing and 
cofinal sequence from $C$. If we write $P(X) := X^p - X - 
b$, we get a sequence $\{a_\rho\}_{\rho < \kappa}$ in $K$ 
such that the sequence $\{v\big(P(a_\rho)\big)\}_{\rho < 
\kappa}$ is stricly increasing and cofinal in $C$. 
Thus, the sequence $\{P(a_\rho)\}_{\rho < \kappa}$ is 
pseudo-convergent (with $0$ one of its limits). As $v
\big(P(a_\alpha)\big) < 0$, we have $v(a_\beta - a_\alpha) 
= \frac{1}{p}v\big(P(a_\alpha)\big) = \gamma_\alpha$ for $
\alpha < \beta < \kappa$. Thus, the sequence $\{a_\rho\}
_{\rho < \kappa}$ is also pseudo-convergent. Furthermore, $\{a_\rho\}_{\rho < \kappa}$ has no limit in $K$: if $l \in K$ is a limit of $\{a_\rho\}_{\rho < \kappa}$ then $P(l)$ is a limit of $\{P(a_\rho)\}_{\rho < \kappa}$. As $\{v\big(P(a_\rho)\big)\}_{\rho < \kappa}$ is cofinal in $C$, $v\big(P(l)\big)$ would be a maximal element of $C$: contradiction.\\[0.1in]
%
Step 3: (Claim) Let $P_0(X) \in K[X]$, and assume that
$v\big(P_0(a_\alpha)\big)$ is strictly increasing ultimately. Then
$\deg\big(P_0(X)\big) \geq p$. \\
We take such a $P_0$ of minimal degree, assume this degree is $n < p$,
and will derive a contradiction. One consequence of Lemma 8 in \cite{Ka} is that: \begin{equation}
 v\big(P_0(a_\rho)\big) = \delta' + \gamma_\rho  \text{ ultimately for } \rho < \kappa
\end{equation} where $\delta'$ is the ultimate valuation of
$P_0'(a_\rho)$ and $\gamma_\rho$ is the valuation of $(a_\sigma -
a_\rho)$ for $\rho < \sigma < \kappa$ (which does not depend on $\sigma$
as $\{a_\rho\}_{\rho < \kappa}$ is pseudo-convergent). We write $P(X) =
\sum_{i=0}^{m} h_i(X)P_0(X)^i $ with $\deg(h_i) < n, \forall i \in \{1,
\dots, m\}$. Then, $\{h_i(a_\rho)\}_{\rho < \kappa}$ is ultimately of
constant valuation, and we let $\lambda_i$ be this valuation. As
$\{a_\rho\}_{\rho < \kappa}$ has no limit in $K$, it is easy to see that
$n>1$, so that $m < p$. By Lemma 4 in \cite{Ka}, 
 there is an integer $i_0 \in \{1, \dots, m \}$ such that we have ultimately: \begin{equation}
\forall i \neq i_0 \quad (\lambda_i + i \delta') + i \gamma_\rho  > (\lambda_{i_0} + i_0 \delta') + i_0 \gamma_\rho.
\end{equation}
Then, ultimately:
\begin{equation}
p\gamma_\rho = v\big(P(a_\rho)\big) = v\Big(\sum_{i=0}^{m} h_i(a_\rho)P_0(a_\rho)^i\Big) = \lambda_{i_0} + i_0(\delta' + \gamma_\rho).
\end{equation}
Thus, we have ultimately $(p - i_0)\gamma_\rho = \lambda_{i_0} + i_0\delta'$. As $p > m \geq i_0$, the left hand side of the equation increases strictly monotonically with $\rho$. But the right hand side is constant: it has no dependence in $\rho$! We have a contradiction, thus $n = p$.\\[0.1in]
Step 4: Clearly, $\{a_\rho\}_{\rho < \kappa}$ is of algebraic type. By
Theorem 3 in \cite{Ka}, if $a_\infty$ is a root of $P$, we get an immediate extension $(L',v') = (K(a_\infty),v')$. Let $a_\infty = a$, we have $(K(a),v')$ isomorphic to $(K(a),\tilde{v})$. Thus:
\begin{equation}
\forall Q \in K_p[X] \quad
\tilde{v}\big(Q(a)\big) = v'\big(Q(a)\big) = v\big(Q(a_\rho)\big) \text{ ultimately}
\end{equation}
This shows the uniqueness of $\tilde{v}$ and concludes the proof of the theorem.
\end{proof}

\bibliographystyle{plain}

\begin{thebibliography}{}

\end{thebibliography}


\begin{thebibliography}{1}




\bibitem{EP}
  A.J. Engler, A. Prestel,
  \emph{Valued fields}.
  Springer,
  2005. 
 

\bibitem{Ka}
 I. Kaplansky,
 \emph{Maximal Fields with Valuation}.
 Duke Math. J. Volume 9, Number 2 (1942), 303 -- 321.
 
\bibitem{Ko}
 J. Koenigsmann,
 \emph{$p$-Henselian Fields}.
 Manuscripta Math. 87 (1995), no. 1, 89 -- 99.
 

\end{thebibliography}

\smallskip\noindent
Adresses and contacts of the authors: \\
DMA, ENS - 45 rue d'Ulm, 75230 Paris cedex 05,
FRANCE\\
e-mail: {\tt zoe.chatzidakis@ens.fr}, {\tt milan.perera@ens.fr}\\

\end{document}